\documentclass[11pt, a4paper]{article}

\usepackage[utf8]{inputenc}
\usepackage{amsmath}
\usepackage{amssymb}
\usepackage{amsthm}
\usepackage{geometry}
\usepackage[colorlinks=true, allcolors=blue]{hyperref}

\newtheorem{theorem}{Theorem}[section]

\newtheorem{definition}[theorem]{Definition}
\newtheorem{proposition}[theorem]{Proposition}

\theoremstyle{remark}
\newtheorem{remark}[theorem]{Remark}
\newtheorem*{theorem*}{Theorem}

\newcommand{\Prob}{\mathbb{P}}
\newcommand{\E}{\mathbb{E}}
\newcommand{\Var}{\mathrm{Var}}

\title{A Probabilistic Framework for the Erdős-Kac Theorem}
\author{Mantha Sai Gopal \\ 
Sri Sathya Sai Institute of Higher Learning, India \\
\texttt{manthasaigopal@gmail.com}}
\date{September 4, 2025}

\begin{document}

\maketitle

\begin{abstract}
    The Erd\H{o}s-Kac theorem, a foundational result in probabilistic number theory, states that the number of prime factors of an integer follows a Gaussian distribution. In this paper we develop and analyze probabilistic models for “random integers” in order to study the mechanisms underlying this theorem. We begin with a simple model, where each prime $p$ is chosen as a divisor with probability $1/p$ in a sequence of independent trials. A preliminary analysis reveals that this construction almost surely yields an integer with infinitely many prime factors. To create a tractable framework, we study a truncated version, $N_x = \prod_{p \le x} p^{X_p}$, where $X_p$ are independent Bernoulli($1/p$) random variables. We prove an analogue of the Erd\H{o}s-Kac theorem within this framework, showing that the number of prime factors $\omega(N_x)$ satisfies a central limit theorem with mean and variance asymptotic to $\log\log x$.
\end{abstract}

\section{Introduction}
The study of the integers, $\mathbb{Z}$, has traditionally been regarded as a deterministic domain. However, in the early twentieth century, a novel perspective emerged. Certain ostensibly irregular features of integers could be elucidated through probabilistic reasoning. This perspective gave rise to probabilistic number theory, a discipline concerned with the statistical behavior of arithmetic functions. Its origins trace back to the seminal work of Hardy and Ramanujan \cite{Hardy1917}, who investigated the function $\omega(n)$, enumerating the number of distinct prime factors of an integer $n$. They established that $\omega(n)$ possesses a normal order of $\log\log n$. Mathematically, for almost all integers $n$, $\omega(n)$ is asymptotically close to $\log\log n$. More precisely, for any $\epsilon > 0$, the set of integers $n$ satisfying

$$
|\omega(n) - \log\log n| > \epsilon \log\log n
$$

has an asymptotic density zero. This result demonstrates that, while $\omega(n)$ may exhibit significant local fluctuations, it conforms to a strikingly regular pattern on the macroscopic scale.

While Hardy and Ramanujan quantified the typical magnitude of $\omega(n)$, the natural subsequent question concerns its distribution. How are the values of $\omega(n)$ distributed around their normal order? The answer to this question was provided by Erdős and Kac in 1940 \cite{Erdos1940}, and constitutes a cornerstone of probabilistic number theory. The Erdős–Kac theorem asserts that, for each real $z$,

$$
\lim_{x \to \infty} \#\Big\{ n \le x : \frac{\omega(n) - \log\log n}{\sqrt{\log\log n}} \le z \Big\} \cdot \frac{1}{x} = \Phi(z),
$$

where $\Phi(z) = \frac{1}{\sqrt{2\pi}} \int_{-\infty}^{z} e^{-t^2/2}\, dt$ denotes the standard normal distribution function.

The theorem suggests that, in a statistical sense, divisibility by distinct primes behaves as though nearly independent. Formally, if one models the event that a large integer is divisible by a prime $p$ as occurring with probability $1/p$, independently across primes, then $\omega(n)$ corresponds to a sum of independent Bernoulli random variables. The classical Central Limit Theorem then predicts a Gaussian distribution for $\omega(n)$, offering a compelling heuristic for the Erdős–Kac phenomenon.

The original proofs employ sophisticated tools such as the Turán–Kubilius inequality \cite{Tenenbaum1995} and sieve methods to control the subtle dependencies imposed by the multiplicative structure of $\mathbb{Z}$. This motivates the question - Can the Gaussian behavior of $\omega(n)$ be derived within a simpler, purely probabilistic framework that ignores these dependencies by design?

In this paper, we construct and analyze such a probabilistic model. Let $\mathcal{P}$ denote the set of all primes, and define a random integer

$$
N = \prod_{p \in \mathcal{P}} p^{X_p},
$$

where the $X_p$ are independent $\mathrm{Bernoulli}(1/p)$ random variables. Our first result exposes a fundamental incompatibility between this independence assumption and the finiteness of integers:

\begin{theorem*}
A random integer $N = \prod_{p \in \mathcal{P}} p^{X_p}$, with the $X_p$ independent $\mathrm{Bernoulli}(1/p)$ variables, possesses infinitely many prime factors almost surely.
\end{theorem*}

This demonstrates that, while the independence assumption facilitates analytic tractability, it fails to enforce the global constraint that integers are finite. To avoid this, we introduce a truncated model, restricting to primes $p \le x$:

$$
N_x = \prod_{p \le x} p^{X_p}.
$$

Within this framework, we establish an analogue of the Erdős–Kac theorem:

\begin{theorem*}[Erdős–Kac Analogue]
Let $N_x = \prod_{p \le x} p^{X_p}$ with $X_p$ independent $\mathrm{Bernoulli}(1/p)$ variables, and define $\Omega_x = \omega(N_x) = \sum_{p \le x} X_p$. Denote $\mu_x = \mathbb{E}[\Omega_x]$ and $\sigma_x^2 = \mathrm{Var}(\Omega_x)$. Then, as $x \to \infty$,

$$
\mu_x = \log\log x + O(1), \quad \sigma_x^2 = \log\log x + O(1),
$$

and the normalized variable converges in distribution to a standard normal

$$
\frac{\Omega_x - \mu_x}{\sigma_x} \xrightarrow{d} \mathcal{N}(0,1).
$$
\end{theorem*}

The paper is organized as follows. Section 2 defines the independent random sieve model and establishes its degeneracy. Section 3 develops the truncated model and proves the Erdős–Kac analogue.

\section{The Random Sieve Model and its Degeneracy}

Let $\mathcal{P} = \{p_1, p_2, \dots\} = \{2, 3, \dots\}$ be the set of prime numbers. We construct a probability space to model the formation of a random square-free integer.

\begin{definition}[The Random Sieve]
Let $(\Omega, \mathcal{F}, \Prob)$ be a probability space. For each prime $p \in \mathcal{P}$, let $X_p$ be an independent Bernoulli random variable such that
$$X_p = \begin{cases} 1 & \text{with probability } 1/p \\ 0 & \text{with probability } 1 - 1/p \end{cases}$$
We define the formal \textbf{random integer} $N$ as the product
$$N = \prod_{p \in \mathcal{P}} p^{X_p}.$$
\end{definition}

The number of distinct prime factors of $N$ is the random variable $\omega(N) = \sum_{p \in \mathcal{P}} X_p$. We now show that this construction almost surely fails to produce a finite integer.

\begin{theorem}\label{thm:infinite}
The random integer $N$ is infinite with probability 1.
\end{theorem}
\begin{proof}
The integer $N$ is finite if and only if the set $\{p \in \mathcal{P} \mid X_p = 1\}$ is finite. This is equivalent to the statement that the sum $\omega(N) = \sum_{p \in \mathcal{P}} X_p$ converges. We appeal to the Borel-Cantelli lemmas. The events $\{X_p = 1\}$ are independent for all $p$. We consider the sum of their probabilities:
$$\sum_{p \in \mathcal{P}} \Prob(X_p = 1) = \sum_{p \in \mathcal{P}} \frac{1}{p}.$$
This is the sum of the reciprocals of the primes, which is known to diverge (cf. Mertens' second theorem \cite{Mertens1874}). By the second Borel-Cantelli lemma, since the events are independent and the sum of their probabilities diverges, infinitely many of the events $\{X_p=1\}$ will occur almost surely. Therefore, $N$ has infinitely many prime factors with probability 1, and thus is infinite.
\end{proof}

\begin{remark}
This result reveals a fundamental limitation of our simple model. The independence assumption, while powerful, fails to enforce the global constraint that integers must be finite.
\end{remark}

\section{The Erd\H{o}s-Kac Theorem for the Independent Model}

To proceed with a well-defined model, we analyze a finite, truncated version that is guaranteed to produce a conventional integer.

\begin{definition}
For any $x > 0$, we define the \textbf{truncated random integer} $N_x$ as
$$N_x = \prod_{p \le x} p^{X_p},$$
where $X_p$ are independent Bernoulli variables. We denote the number of prime factors of $N_x$ by the random variable
$$\Omega_x = \omega(N_x) = \sum_{p \le x} X_p.$$
\end{definition}

$\Omega_x$ is a sum of a finite number of independent (but not identically distributed) Bernoulli random variables. We begin by computing its first two moments, which mirror the classical results for $\omega(n)$.

\begin{proposition}[Moments of $\Omega_x$] \label{thm:moments}
Let $\mu_x = \E[\Omega_x]$ and $\sigma_x^2 = \Var(\Omega_x)$. Then as $x \to \infty$,
\begin{align*}
    \mu_x &= \log\log x + O(1), \\
    \sigma_x^2 &= \log\log x + O(1).
\end{align*}
\end{proposition}

\begin{proof}
By the linearity of expectation,
$$\mu_x = \E\left[\sum_{p \le x} X_p\right] = \sum_{p \le x} \E[X_p] = \sum_{p \le x} \frac{1}{p}.$$
By Mertens' second theorem, this is $\log\log x + B_1 + o(1)$, where $B_1$ is the Meissel-Mertens constant.

For the variance, since the $X_p$ are independent, the variance of the sum is the sum of the variances:
$$\sigma_x^2 = \Var\left(\sum_{p \le x} X_p\right) = \sum_{p \le x} \Var(X_p).$$
For a Bernoulli variable with success probability $q$, the variance is $q(1-q)$. Thus,
$$\sigma_x^2 = \sum_{p \le x} \frac{1}{p}\left(1 - \frac{1}{p}\right) = \sum_{p \le x} \frac{1}{p} - \sum_{p \le x} \frac{1}{p^2}.$$
The first term is $\log\log x + O(1)$. The second sum, $\sum_{p \le x} 1/p^2$, converges to the prime zeta function value $P(2) \approx 0.4522$ as $x \to \infty$. Therefore,
$$\sigma_x^2 = (\log\log x + O(1)) - P(2) = \log\log x + O(1).$$
The proposition is proved.
\end{proof}

We now prove our main result for this model, an analogue of the Erd\H{o}s-Kac theorem.

\begin{theorem}[Erd\H{o}s-Kac Analogue]\label{thm:clt-independent}
The random variable $\Omega_x = \sum_{p \le x} X_p$, when normalized, converges in distribution to a standard normal random variable. That is,
$$\frac{\Omega_x - \mu_x}{\sigma_x} \xrightarrow{d} \mathcal{N}(0, 1) \quad \text{as } x \to \infty,$$
where $\mu_x \sim \log\log x$ and $\sigma_x^2 \sim \log\log x$.
\end{theorem}

\begin{proof}
We have a sum of independent, non-identically distributed random variables, thus we shall use the Lindeberg-Feller Central Limit Theorem \cite{Feller1971}. Let $Y_p = X_p - \E[X_p] = X_p - 1/p$. Then $\Omega_x - \mu_x = \sum_{p \le x} Y_p$. The variance is $\sigma_x^2 = \sum_{p \le x} \E[Y_p^2]$. The Lindeberg condition states that for any $\epsilon > 0$,
$$\lim_{x \to \infty} \frac{1}{\sigma_x^2} \sum_{p \le x} \E[Y_p^2 \cdot \mathbf{1}_{\{|Y_p| > \epsilon \sigma_x\}}] = 0,$$
where $\mathbf{1}_{\{\cdot\}}$ is the indicator function.

Let us analyze the terms $Y_p$. The random variable $Y_p$ takes the value $1 - 1/p$ with probability $1/p$, and the value $-1/p$ with probability $1 - 1/p$. In either case, $|Y_p| \le 1$.
From Proposition \ref{thm:moments}, the variance $\sigma_x^2 \sim \log\log x$, which tends to infinity as $x \to \infty$. Thus, for any fixed $\epsilon > 0$, we can find a $X_0$ such that for all $x > X_0$, we have $\epsilon \sigma_x > 1$.
For such $x$, the condition $|Y_p| > \epsilon \sigma_x$ can never be met, since $|Y_p| \le 1$. The indicator function $\mathbf{1}_{\{|Y_p| > \epsilon \sigma_x\}}$ is therefore 0 for all $p \le x$.

Consequently, for any $x > X_0$, the sum in the Lindeberg condition is exactly 0:
$$\sum_{p \le x} \E[Y_p^2 \cdot \mathbf{1}_{\{|Y_p| > \epsilon \sigma_x\}}] = \sum_{p \le x} \E[Y_p^2 \cdot 0] = 0.$$
The limit is therefore 0, and the Lindeberg condition is satisfied. By the Lindeberg-Feller CLT, the normalized sum converges in distribution to $\mathcal{N}(0,1)$.
\end{proof}

\end{document}